\documentclass{amsart}

\input xy
\xyoption {all}

\usepackage{latexsym} 
\usepackage{amssymb} 
\usepackage{color} 
\usepackage[frame,cmtip,arrow,matrix,line,graph]{xy}

\usepackage{amsmath} 
\usepackage{amsthm} 
\usepackage{amsfonts} 
\usepackage{amscd}
 
\usepackage{epsfig} 

\usepackage{caption2}
\usepackage{mathrsfs}

\newtheorem{thm}{Theorem}

\newcommand{\Dif}{{\rm{Diff}} }

\newcommand{\Ñ}{\kern-.1pt\vrule height2.8ptwidth3ptdepth-2.3pt\kern1pt}

\newcommand{\Sph}{\mathcal{S}}

\newcommand{\del}{\partial}

\newcommand{\cln}{\colon\!}

\newcommand{\lgl}{\langle}
\newcommand{\rgl}{\rangle}

\newcommand{\x}{\!\times\!}

\newcounter{samcounter}

\begin{document}

\title{Erratum to ``Stabilization for the automorphisms of free groups with boundaries''}

\author{Allen Hatcher}
\author{Nathalie Wahl}

\maketitle

The purpose of this note is to point out a gap in an argument in our paper \cite{HW} and explain how to fill it (or bypass it).

The gap in \cite{HW} occurs in Section~4 in the proof of assertions (A) and (B),
at the point where there are diagram chasing arguments in two commutative diagrams (displayed on page 1333 for case (A)). In each diagram the groups $G_n$ in the two rows are isomorphic but not identical. 
If we denote by $G_n$ and $G_n'$ these two isomorphic groups, the first diagram chase needs the composition $H_i(G_n)\to H_i(G_{n+1})\to H_i(G_{n+1},G_n')$ to be trivial, which is the case because $G_n$ and $G_n'$ are conjugate in $G_{n+1}$. The second diagram chase needs the composition 
$H_{i+1}(G_n,G_{n-1})\to H_i(G_{n-1})\to H_i(G_n')$ to be trivial, but there is no {\it a priori\/} reason for this to be true, although it is true and follows {\it a posteriori\/} from Theorem~\ref{alpha} below.
An analogous diagram chase is used in the proof of assertion (C), but in that case there is an isomorphism $G_n\to G_n'$ commuting with the inclusion of $G_{n-1}$, so that the diagram chase is correct. For (A) and (B), there is an isomorphism $G_n\to G_n'$ that commutes with the inclusion of a subgroup $H$ of $G_{n-1}$, with $H_*(H)\cong H_*(G_{n-1})$ in a range given by Theorem~\ref{alpha}.

Recall from \cite{HW} that
$M_{n,k}^s=N\# (\#_n S^1\x S^2)\# (\#_k S^1\x D^2)\# (\#_s D^3)$, where $N$ is a fixed compact connected oriented 3-manifold,
and that $A_{n,k}^s$ denotes the quotient of the mapping class group $\pi_0\Dif(M_{n,k}^s\,\textrm{rel}\,\del M_{n,k}^s)$ by twists along spheres embedded in $M$. 
Assertion (A) says that the map $\alpha_i\cln H_i(A_{n,k}^{s+2})\to H_i(A_{n+1,k}^{s+1})$, induced by identifying discs in the last two boundary spheres of $M$, 
is surjective when $n\ge 3i$ and an isomorphism when $n\ge 3i+2$. Assertion (B) says that the same is true for the map  
$\beta_i\cln H_i(A_{n,k}^{s+2})\to H_i(A_{n+1,k}^{s})$, where $\beta_i$ is induced by identifying the last two boundary spheres of $M$. 
Assertions (A) and (B) are used in \cite{HW} to show (1) and (3) in Theorem 4.1, namely that the stabilization maps $\phi_i\cln H_i(A_{n,k}^{s})\to H_i(A_{n+1,k}^{s})$, induced by gluing a twice punctured $S^1\times S^2$ along a boundary sphere, and $\mu_i\cln H_i(A_{n,k}^{s})\to H_i(A_{n,k}^{s+1})$, induced by gluing a 3-punctured sphere along one of its boundary components, are isomorphisms when $n\ge 3i+3$. To fix the gap described above, it is sufficient to know --- independently --- that $\mu_i\cln H_i(A_{n-1,k}^{s+1})\to H_i(A_{n-1,k}^{s+2})$ is surjective  when $n-1\ge 3i+1$, but the new argument actually proves simultaneously that $\phi_i$ and $\mu_i$ are isomorphisms in the above range,  
 avoiding completely the two commutative diagrams mentioned above, so that the point becomes moot. 
The proof, sketched below, uses two simplicial complexes built from a sphere complex studied in \cite{HW}. The main ingredient of the proof is the high connectivity of these complexes, which is deduced from the connectivity of the earlier sphere complex by a combinatorial argument. This last argument is applied in a more general context in \cite{HW2} and hence not repeated here. 
We actually obtain in this way a better stability range, namely an isomorphism when $n\ge 2i+2$ for both $\phi_i$ and $\mu_i$.

\medskip

Let $M=M_{n,k}^s$ be as above and let $x_0,x_1$ be two points in two boundary spheres $\del_0M,\,\del_1M$ of $M$, where we allow the possibility $x_0=x_1$ with $\del_0M=\del_1M$. Recall from \cite{HW} that $\Sph_c(M)$ denotes the simplicial complex of isotopy classes of non-separating sphere systems in $M$. We use here an enhanced version of $\Sph_c(M)$: Let $X^A(M,x_0,x_1)$ be the simplicial complex whose vertices are pairs $(S,a)$, where $S$ is in $\Sph_c(M)$ and $a$ is the isotopy class of an arc from $x_0$ to $x_1$ intersecting $S$ transversely in exactly one point, with a choice of orientation of the arc if $x_0=x_1$. A $p$-simplex of $X^A(M,x_0,x_1)$ is a collection $\lgl (S_0,a_0),\dots,(S_p,a_p)\rgl$ of such pairs that are disjoint except for the endpoints of the arcs, and such that the spheres $\lgl S_0,\dots,S_p\rgl$ form a $p$-simplex of $\Sph_c(M)$.

There is a map $X^A(M,x_0,x_1)\to \Sph_c(M)$ which forgets the arcs, and we think of $X^A(M,x_0,x_1)$ as the complex $\Sph_c(M)$ with arcs labeling its vertices. In \cite{HW2}, we show how the connectivity of $X^A(M,x_0,x_1)$ can be deduced from that of $\Sph_c(M)$, established in \cite{HW} (Proposition~3.2), using a combinatorial argument \cite[Thm.~3.8]{HW2}. 

\begin{thm}\cite[Prop.~4.4]{HW2}
$X^A(M_{n,k}^s,x_0,x_1)$ is $(\frac{n-3}{2})$-connected.
\end{thm}

The group $A_{n,k}^s$ acts on the complex $X^A(M_{n,k}^s,x_0,x_1)$ since twists along spheres act trivially on embedded spheres and arcs that meet these spheres in one point transversely. 
The stabilizer of a vertex is isomorphic to $A_{n-1,k}^s$ both in the case that $x_0=x_1$ and in the case that $x_0$ and $x_1$ lie on different boundary components of $M$, but in the first case the inclusion of the stabilizer into $A_{n,k}^s$ induces the stabilization $\phi_i$ described above, while in the second case it is the map $\psi\cln A_{n-1,k}^s\to A_{n,k}^s$ induced by gluing a 4-punctured sphere along two of its boundary spheres. The spectral sequences for  the action of $A_{n,k}^s$ on $X^A(M_{n,k}^s,x_0,x_1)$ in each of the two cases yield the following stability:

\begin{thm}\cite[Thm.~6.1]{HW2}\label{alpha}
The maps $\phi_i,\psi_i\cln H_i(A_{n,k}^s) \to H_i(A_{n+1,k}^s)$ are isomorphisms when $n\ge 2i+2$ and surjections when $n\ge 2i+1$ {\rm (}with $s\ge 1$ for $\phi$ and $s\ge 2$ for $\psi${\rm \/)}. 
\end{thm}
(In \cite{HW2}, the two maps are denoted $\alpha$ and $\beta$.)
Note now that $\phi=\eta\mu$ and $\psi=\mu\eta$, for $\eta$ the map obtained by gluing a 3-punctured sphere along two of its boundary spheres. It follows from the theorem that both $\mu$ and $\eta$ are isomorphisms in the same range. To show that $\mu$ is surjective in that range, which is enough to fix the gap, it suffices to consider only the map $\psi$, that is, the case where $x_0$ and $x_1$ lie on two different boundary components.

\medskip

An analogous gap occurred in \cite{HV3}, in the proof of homological stability for $\operatorname{Out}(F_n)$, which is the group $A_{n,0}^0$ with $N=S^3$ in our notation. The paper \cite{HVW} bypasses the gap by showing that the complex $X^A(M_{n,0}^s,x_0,x_1)$ (denoted $Z_n$ in that paper) with $x_0$ and $x_1$ on different boundary components of $M_{n,0}^s$, becomes contractible when $n$ goes to infinity, which is enough in that case by the main result of \cite{HV1}. The present erratum gives an alternative argument to fill in the gap of \cite{HV3}, which does not depend on \cite{HV1}.

\end{document}